\documentclass[a4paper]{amsart}
\usepackage{amsmath,amsthm,amssymb,latexsym,epic,bbm,comment,mathrsfs}
\usepackage{graphics,enumerate,stmaryrd,color,tikz,ytableau}
\usepackage[all,2cell]{xy}
\xyoption{2cell}
\usepackage[most]{tcolorbox}

\newtheorem{theorem}{Theorem}
\newtheorem{lemma}[theorem]{Lemma}
\newtheorem{corollary}[theorem]{Corollary}
\newtheorem{proposition}[theorem]{Proposition}

\usepackage[all]{xy}
\usepackage[active]{srcltx}
\usepackage[parfill]{parskip}
\usepackage{enumerate}

\def\HK#1{\textcolor[rgb]{.75,0.00,0.00}{[HK: #1]}}

\begin{document}

\title[Addendum to: Some homological 
properties of category $\mathcal{O}$, V]
{Addendum and Corrigendum to: \\Some homological 
properties of category $\mathcal{O}$, V}

\author[Hankyung Ko and Volodymyr Mazorchuk]
{Hankyung Ko and Volodymyr Mazorchuk}

\begin{abstract}
In various special cases, we provide an alternative argument for 
equivalence between Kostant positivity of some involutions in type $A$
and the conjunction of the categorical K{\aa}hrstr{\"o}m's condition
for these involutions with indecomposability of the images of the 
corresponding simple modules in category $\mathcal{O}$ under projective functors.
\end{abstract}

\maketitle


\section{Introduction and description of the results}\label{s1}

\subsection{Introduction}\label{s1.1}

The elements of the symmetric group $S_n$ index bijectively, $w\mapsto L_w$, 
the isomorphism classes of simple objects in 
the principal block of the BGG category $\mathcal{O}$ 
for $\mathfrak{sl}_n$. For each such $L_w$, one can ask a question
called {\em Kostant's problem}, see Section~\ref{s3.1} and \cite{Jo}, about the adjoint action of 
$\mathfrak{sl}_n$ on the associative algebra of linear endomorphisms of
$L_w$. The answer can be either positive or negative and we accordingly call
both $L_w$ and 
$w$ {\em Kostant positive} or {\em negative}. So far, there does not seem to be
any reasonable combinatorial criterion for determination of this answer for
a given $w$.

In \cite{KMM}, one can find a conjecture, due to Johan K{\aa}hrstr{\"o}m,
which, in the case when $w$ is an involution, relates the answer to 
Kostant's problem for $L_w$ to various categorical properties of 
the images of $L_w$ under projective functors and also to  combinatorial 
properties of the (dual) Kazhdan-Lusztig basis elements of the Hecke algebra.
Moreover, \cite[Theorem~8.16]{KMM} asserts that Kostant positivity of 
an involution $w\in S_n$ is equivalent to the conjunction of the 
categorical K{\aa}hrstr{\"o}m's condition for $w$ with indecomposability 
of the images of $L_w$ under projective functors. However, 
Eduardo Venturini has recently pointed out a gap in the proof of 
\cite[Proposition~8.15]{KMM}, which is the key statement behind 
one of the implications in \cite[Theorem~8.16]{KMM}.
The problem is that the functor $\Theta$ defined in the proof 
sends $\theta_xL_d$ to $\theta_x \otimes_U \mathrm{Res}^{A_d}_{U_d}(A_d)$
and not to $\theta_x \otimes_U U_d$ as claimed.

At the moment, we do not see how to fill this gap in the general case.
However, we can fill it in many special cases, in particular, in all
special cases in which this result was used in the literature so far,
that is, in the setups of \cite{KMM,MMM,CM2,CM4}.

\subsection{Results}\label{s1.2}

Our main result is the following special case of \cite[Theorem~8.16]{KMM}. 
See Corollary~\ref{cor:kmmfc}, Proposition~\ref{prop:shape1},
Proposition~\ref{prop:shape2} and all subsections in Section~\ref{s5}
which treat the small rank cases one by one.

{\bf Theorem A.}
{\em
Let $d\in S_n$ be an involution such that one of the following
conditions is satisfied:

\begin{enumerate}[$($a$)$]
\item\label{ca} $d$ is fully commutative;
\item\label{cc} the Robinson-Schensted correspondent of $d$
has shape $(n-2,1^2)$ or $(n-3,2,1)$;
\item\label{cb} $n\leq 7$.
\end{enumerate}
Then the module $L_d$ is Kostant positive if and only if the 
following two conditions are satisfied:
\begin{enumerate}[$($i$)$]
\item\label{ci} applying to $L_d$ any indecomposable projective 
endofunctor of the principal block of category $\mathcal{O}$
outputs either an indecomposable module or zero;
\item\label{cii} if $\theta_x$ and $\theta_y$ are two non-isomorphic 
indecomposable projective 
endofunctors of the principal block of category $\mathcal{O}$
such that both $\theta_x L_d\neq 0$ and $\theta_y L_d\neq 0$,
then $\theta_x L_d\not\cong \theta_y L_d$.
\end{enumerate}
}

We note that Condition~\eqref{ca} covers  the situation
considered in \cite{MMM},
Condition~\eqref{cb} covers  the situation
considered in \cite{KMM,CM4}, and
Condition~\eqref{cc} covers  the situation
considered in \cite{CM2}.

That the failure of either Condition~\eqref{ci} or Condition~\eqref{cii}
forces $L_d$ to be Kostant negative is true for any $d\in S_n$,
see \cite[Proposition~8.14]{KMM}.

In the case of fully commutative involutions, our argument
is quite general. For all other special cases, our argument
is based on a case-by-case work-through, which turns out
to be quite  substantial for shape $(n-3,2,1)$.

\subsection*{Acknowledgments}

We thank Eduardo Venturini for pointing out the abovementioned gap in the 
proof of \cite[Proposition~8.15]{KMM}.

\section{Setup}\label{s2}

\subsection{Category $\mathcal{O}$}\label{s2.1}

Let $\mathfrak{g}=\mathfrak{sl}_n(\mathbb{C})$
with a fixed standard triangular decomposition
\begin{displaymath}
\mathfrak{g}=\mathfrak{n}_-\oplus \mathfrak{h}\oplus\mathfrak{n}_+. 
\end{displaymath}
Let $\mathcal{O}$ be the corresponding BGG category $\mathcal{O}$
and $\mathcal{O}_0$ the principal block of $\mathcal{O}$,
see \cite{BGG,Hu}.

Let $\mathbf{R}$ be the root system of $(\mathfrak{g},\mathfrak{h})$.
Let $W$ be the Weyl group of $(\mathfrak{g},\mathfrak{h})$
which acts on $\mathfrak{h}^*$ in the usual way and via the 
dot-action $w\cdot \lambda:=w(\lambda+\rho)-\rho$, where 
$\rho$ is the half-sum of positive roots.
We denote the Bruhat order on $W$ by $\preceq$.

For $w\in W$, we denote by:
\begin{itemize}
\item $L_w:=L(w\cdot 0)$ the simple highest weight module with 
highest weight $w\cdot 0$;
\item $\Delta_w$ the Verma cover of $L_w$ in $\mathcal{O}$;
\item $P_w$ the (indecomposable) projective cover of $L_w$ in $\mathcal{O}$;
\item $I_w$ the injective envelope of $L_w$ in $\mathcal{O}$;
\item $\nabla_w$ the dual Verma envelope of $L_w$ in $\mathcal{O}$;
\item $T_w$ the tilting envelope of $\Delta_w$ in $\mathcal{O}$.
\end{itemize}
We denote by $\star$ the usual simple preserving duality on
$\mathcal{O}$. In particular, we have $L_w^\star\cong L_w$,
$T_w^\star\cong T_w$, $\Delta_w^\star\cong \nabla_w$
and $P_w^\star\cong I_w$.

\subsection{Projective functors}\label{s2.2}

Let $\mathscr{P}$ be the monoidal category of projective
endofunctors of $\mathcal{O}_0$, see \cite{BG}.
For $w\in W$ we denote by $\theta_w\in \mathscr{P}$
the unique (up to isomorphism) indecomposable object
such that $\theta_w P_e\cong P_w$.

\subsection{Kazhdan-Lusztig combinatorics}\label{s2.3}

Let $\mathbf{H}$ be the Hecke algebra of $W$ with coefficients
in $
\mathbb{Z}[v,v^{-1}]$. Let
\begin{itemize}
\item $\{H_w\,:\,w\in W\}$ be the standard basis of $\mathbf{H}$;
\item $\{\underline{H}_w\,:\,w\in W\}$ be the Kazhdan-Lusztig 
basis of $\mathbf{H}$;
\item $\{\underline{\hat{H}}_w\,:\,w\in W\}$ be the dual
Kazhdan-Lusztig  basis of $\mathbf{H}$.
\end{itemize}
Let $p_{x,y}\in \mathbb{Z}[v,v^{-1}]$,
for $x,y\in W$, be the Kazhdan-Lusztig
polynomials defined via:
\begin{displaymath}
\underline{H}_x=\sum_y  p_{x,y} H_x. 
\end{displaymath}
Note that  $p_{x,y}=0$ unless $x\preceq y$.
We refer to \cite{KL,So2} for details. We use the normalization of \cite{So2}.

For $x,y,z\in W$, define the structure constant 
$h_{x,y}^z\in \mathbb{Z}[v,v^{-1}]$ via:
\begin{displaymath}
\underline{H}_x\underline{H}_y=
\sum_{z}h_{x,y}^z\underline{H}_z.
\end{displaymath}
Note that $h_{x,y}^z$ is invariant under the flip $v\mapsto v^{-1}$.
We have Lusztig's $\mathbf{a}$-function 
$\mathbf{a}:W\to\mathbb{Z}_{\geq 0}$, which assigns to
$w\in W$ the maximal $v$-degree of $h_{x,y}^w$, taken 
over all $x,y$. The function $\mathbf{a}$ is constant on 
two-sided cells in $W$. If $w_0^P$ is the longest element of
some Young (parabolic) subgroup $W^P$ of $W$, then 
$\mathbf{a}(w)=\ell(w)$, where $\ell(w)$ is the usual 
Coxeter length of $w$. Note that $\mathbf{a}(w)<\ell(w)$,
for all other $w$.

We also have the Kazhdan-Lusztig $\mu$-function
which is defined  via $\mu(x,y)=0$ if
$x$ and $y$ are not Bruhat comparable and
$\mu(y,x)=\mu(x,y)$ otherwise, with 
$\mu(x,y)$ defined as  the coefficient at $v$
in $p_{x,y}$, if $x\preceq y$.

\subsection{Graded lifts}\label{s2.4}

Let $A$ be a finite dimensional associative algebra such that
$A$-mod is equivalent to $\mathcal{O}_0$. Then $A$ is Koszul
by \cite{So}, in particular, it admits a $\mathbb{Z}$-grading.
We denote by ${}^{\mathbb{Z}}\mathcal{O}_0$ the category of 
finite dimensional $\mathbb{Z}$-graded $A$-modules
(with morphisms being homogeneous maps of degree $0$).
We refer to ${}^{\mathbb{Z}}\mathcal{O}_0$ as the 
graded lift of $\mathcal{O}_0$.
We denote by ${}^{\mathbb{Z}}\mathscr{P}_0$ the 
graded lift of $\mathscr{P}_0$.

We denote by $\langle 1\rangle$ the auto-equivalence of
${}^{\mathbb{Z}}\mathcal{O}_0$ which shifts the grading 
by decreasing the degree by $1$.
For $w\in W$, we fix the following standard graded
lifts for structural objects:
\begin{itemize}
\item $L_w$ the simple object concentrated in degree $0$;
\item $\Delta_w$ has top in degree $0$;
\item $P_w$ has top in degree $0$;
\item $I_w$ has socle in degree $0$;
\item $\nabla_w$ has socle in degree $0$;
\item $T_w$ has the unique $L_w$ subquotient in degree $0$.
\end{itemize}
We denote by $\star$ the standard graded lift of the simple preserving
duality (note that it maps degree $i$ to degree $-i$).

We also denote by $\theta_w$ the graded lift of 
the corresponding indecomposable projective functor
such that $\theta_w P_e\cong P_w$ holds in 
${}^{\mathbb{Z}}\mathcal{O}_0$.

With these conventions, we have the following decategorification result:
we fix an isomorphism between the split Grothendieck ring
$\mathbf{Gr}_\oplus({}^{\mathbb{Z}}\mathscr{P}_0)$
and $\mathbf{H}$ which sends $\theta_w$ to 
$\underline{H}_w$, for $w\in W$, and $\langle 1\rangle$
to $v^{-1}$; 
we fix an isomorphism between the split Grothendieck group
$\mathbf{Gr}_\oplus({}^{\mathbb{Z}}\mathcal{O}_0)$
and $\mathbf{H}$ which sends $\Delta_w$ to 
${H}_w$, for $w\in W$, and $\langle 1\rangle$
to $v^{-1}$. Then the action of ${}^{\mathbb{Z}}\mathscr{P}_0$
on ${}^{\mathbb{Z}}\mathcal{O}_0$ decategorifies to the right regular
action of $\mathbf{H}$. Under this map, $P_w$ gets sent to
$\underline{H}_w$ while $L_w$ gets sent to
$\underline{\hat{H}}_w$.

As usual, we denote the Hom space in ${}^\mathbb{Z}\mathcal{O}_0$
by $\mathrm{hom}$.

\subsection{Cells}\label{s2.5}

For $x,y\in W$, we write $x\geq_L y$ provided that there is some $w\in W$
such that $\underline{H}_x$ appears with a non-zero coefficient when
expressing $\underline{H}_w\underline{H}_y$ in the KL basis.
We set $x\sim_L y$ if and only if $x\geq_L y$ and $y\geq_L x$.
The pre-order $\geq_L$ is called the left Kazhdan-Lusztig pre-order
and the equivalence classes the left Kazhdan-Lusztig cells.
The right versions $\geq_R$ and $\sim_R$ as well as  the
two-sided versions $\geq_J$ and $\sim_J$ are defined similarly
using multiplication on the right and on both sides, respectively.

In our case, we have $W\cong S_n$, the symmetric group on $n$
elements. The Robinson-Schensted correspondence is a bijection
\begin{displaymath}
\mathbf{RS}:S_n\to
\coprod_{\lambda\vdash n}\mathbf{STY}_\lambda\times \mathbf{STY}_\lambda,
\end{displaymath}
given by Schensted's insertion algorithm, see \cite{Sa,Sch}, between 
$S_n$ and pairs of standard Young tableaux of the same shape
(this shape is a partition of $n$). The bijection is written
$w\mapsto (\mathtt{P}_w,\mathtt{Q}_w)$, where $\mathtt{P}_w$
is the insertion tableau and $\mathtt{Q}_w$
is the recording tableau. For $w\in S_n$, we denote by 
$\mathbf{sh}(w)$ the shape of $\mathtt{P}_w$
(equivalently, the shape of $\mathtt{Q}_w$).

By \cite{KL}, for $x,y\in W$, we have:
\begin{itemize}
\item $x\sim_L y$ if and only if $\mathtt{Q}_x=\mathtt{Q}_y$;
\item $x\sim_R y$ if and only if $\mathtt{P}_x=\mathtt{P}_y$;
\item $x\sim_J y$ if and only if $\mathbf{sh}(x)=\mathbf{sh}(y)$.
\end{itemize}
By \cite{Ge}, the order  $\geq_J$ is the opposite of the
dominance order on partitions.

Each left and each right cell contains a unique involution
since $w\in S_n$ satisfies $w^2=e$ if and only if 
$\mathtt{P}_w=\mathtt{Q}_w$.

\subsection{Categories $\mathcal{O}^{\hat{\mathcal{R}}}$}\label{s2.6}

Let $\mathcal{R}$ be a fixed right cell of $W$.
Denote by $\hat{\mathcal{R}}$ the set of all 
elements $w\in W$ satisfying $w\leq_R \mathcal{R}$.
Let $\mathcal{O}^{\hat{\mathcal{R}}}_0$ be the Serre subcategory
of $\mathcal{O}_0$ generated by all $L_w$ where 
$w\in \hat{\mathcal{R}}$.

The exact embedding $\mathcal{O}^{\hat{\mathcal{R}}}_0\hookrightarrow
\mathcal{O}_0$ has both a left adjoint, which we denote by
$\mathrm{Z}^{\hat{\mathcal{R}}}$, and a right adjoint, 
which we denote by $\mathrm{Z}_{\hat{\mathcal{R}}}$.
For $w\in W$, we set 
\begin{displaymath}
P^{\hat{\mathcal{R}}}_w:=\mathrm{Z}^{\hat{\mathcal{R}}}(P_w),\quad
\Delta^{\hat{\mathcal{R}}}_w:=\mathrm{Z}^{\hat{\mathcal{R}}}(\Delta_w),\quad
I^{\hat{\mathcal{R}}}_w:=\mathrm{Z}_{\hat{\mathcal{R}}}(I_w),\quad
\nabla^{\hat{\mathcal{R}}}_w:=\mathrm{Z}_{\hat{\mathcal{R}}}(\nabla_w).
\end{displaymath}
Directly from the definitions, all these modules admit graded lifts
which we denote by the same symbols.
Then $P^{\hat{\mathcal{R}}}_w$ is the projective cover of 
$L_w$  in $\mathcal{O}^{\hat{\mathcal{R}}}_0$
and  $I^{\hat{\mathcal{R}}}_w$ the injective envelope of 
$L_w$ in $\mathcal{O}^{\hat{\mathcal{R}}}_0$. 
Note that $\mathcal{O}^{\hat{\mathcal{R}}}_0$ inherits from 
$\mathcal{O}_0$ a graded lift, which we denote by 
${}^\mathbb{Z}\mathcal{O}^{\hat{\mathcal{R}}}_0$.
Also, $\mathcal{O}^{\hat{\mathcal{R}}}_0$ is stable under the action of
$\mathscr{P}_0$, see \cite[Lemma~13]{MM1}.

For $w\in \hat{\mathcal{R}}$, we have an ungraded isomorphism
$P^{\hat{\mathcal{R}}}_w\cong I^{\hat{\mathcal{R}}}_w$
if and only if $w\in \mathcal{R}$.
We refer to \cite{MS,Ma2} for details.

\section{Harish Chandra bimodules and Kostant's problem}\label{s3}

\subsection{Original Kostant's problem}\label{s3.1}

Let $U(\mathfrak{g})$ be the universal enveloping algebra of 
$\mathfrak{g}$. Let $M$ and $N$ be $\mathfrak{g}$-modules. Then
the vector space $\mathrm{Hom}_\mathbb{C}(M,N)$ has the
structure of $U(\mathfrak{g})$-$U(\mathfrak{g})$-bimodule via the $U(\mathfrak{g})$-module structures on $M,N$.
Let us denote by $\mathcal{L}(M,N)$ the set of elements of 
$\mathrm{Hom}_\mathbb{C}(M,N)$ on which the adjoint action of 
$U(\mathfrak{g})$ is finite. Then
$\mathcal{L}(M,N)$ is a 
$U(\mathfrak{g})$-$U(\mathfrak{g})$-subbimodule
of $\mathrm{Hom}_\mathbb{C}(M,N)$.

The representation map for $M$
embeds $U(\mathfrak{g})/\mathrm{Ann}_{U(\mathfrak{g})}(M)$ into 
$\mathcal{L}(M,M)$. The original Kostant problem
for $M$, as formulated in \cite{Jo}, asks whether this 
embedding of $U(\mathfrak{g})/\mathrm{Ann}_{U(\mathfrak{g})}(M)$ into 
$\mathcal{L}(M,M)$ is an isomorphism.
We write $\mathbf{K}(M)=\mathtt{true}$ 
and say that $M$ is {\em Kostant positive} if the answer is ``yes''.
We write $\mathbf{K}(M)=\mathtt{false}$ 
and say that $M$ is {\em Kostant negative} if the answer is ``no''.

In general, the answer to this problem is not known even for
the modules $L_w$. However, it is known that 
$\mathbf{K}(L_x)=\mathbf{K}(L_y)$ provided that 
$x\sim_L y$ (in our setup of type $A$), see \cite{MS}.
For simple modules in $\mathcal{O}_0$ this, in particular,
reduces the study of Kostant's problem to the modules $L_d$, where 
$d\in W$ is an involution.

\subsection{Reformulation using $\mathcal{O}^{\hat{\mathcal{R}}}$}\label{s3.2}

Let  $d\in W$ be an involution and $\mathcal{R}$ the right cell 
containing $d$. Consider the corresponding category 
$\mathcal{O}^{\hat{\mathcal{R}}}$ and the dominant projective
module $P^{\hat{\mathcal{R}}}_e$. The module 
$P^{\hat{\mathcal{R}}}_e$ has simple socle $L_d$
and any composition subquotient of 
$P^{\hat{\mathcal{R}}}_e/L_d$ is of the form
$L_w$, for $w\in \hat{\mathcal{R}}\setminus \mathcal{R}$,
see \cite{Ma2}. Consequently, $P^{\hat{\mathcal{R}}}_e$
is a submodule of $I^{\hat{\mathcal{R}}}_d$.

Denote by $Q_d$ the maximal submodule of 
$I^{\hat{\mathcal{R}}}_d$
with the property that any composition subquotient of 
$Q_d/L_d$ is of the form
$L_w$, for $w\in \hat{\mathcal{R}}\setminus \mathcal{R}$.
We grade $Q_d$ so that it is a submodule of $I^{\hat{\mathcal{R}}}_d \langle -\mathbf{a}(d)\rangle$.
Then we have a graded inclusion $P^{\hat{\mathcal{R}}}_e\subset Q_d$.
By \cite{KhM}, we have $\mathbf{K}(L_d)=\mathtt{true}$
if and only if $Q_d= P^{\hat{\mathcal{R}}}_e$.

\subsection{Harish Chandra bimodules}\label{s3.3}

Let $Z(\mathfrak{g})$ be the center of $U(\mathfrak{g})$.
We denote by $\chi:Z(\mathfrak{g})\to\mathbb{C}$ the central
character of the trivial $\mathfrak{g}$-module $L_e$.
Let $\mathbf{m}$ be the kernel of $\chi$.

Let $\mathcal{H}$ denote the category of Harish Chandra
$U(\mathfrak{g})$-$U(\mathfrak{g})$-bimodules, that is,
finitely generated bimodules the adjoint action of 
$U(\mathfrak{g})$ on which is locally finite and has finite
multiplicities, see \cite{Ja}. Let
${}^\infty_{\,\,0} \mathcal{H}_0^1$ be the full subcategory of
$\mathcal{H}$ which consists of all objects the right action 
of $\mathbf{m}$ on which is zero and the left action of 
$\mathbf{m}$ on which is locally nilpotent. Then
${}^\infty_{\,\,0} \mathcal{H}_0^1$ is equivalent to
$\mathcal{O}_0$ via the pair of adjoint equivalences
$B\mapsto B\otimes_{U(\mathfrak{g})}\Delta_e$
and $M\mapsto \mathcal{L}(\Delta_e,M)$, see \cite{BG}.
For $w\in S_n$, the indecomposable projective functor $\theta_w$ 
corresponds to an indecomposable projective object in
${}^\infty_{\,\,0} \mathcal{H}_0^1$ under this equivalence.
Therefore we will abuse notation and denote the 
indecomposable projective objects in
${}^\infty_{\,\,0} \mathcal{H}_0^1$ by $\theta_w$.
Note that $\theta_w$ corresponds to $P_w$ via $\theta_w P_e\cong P_w$.
Also, due to our equivalence, we can consider a graded 
version of ${}^\infty_{\,\,0} \mathcal{H}_0^1$ which we will denote
in the same way.

We have the full subcategory 
${}^1_{0} \mathcal{H}_0^1$ of 
${}^\infty_{\,\,0} \mathcal{H}_0^1$
consisting of all objects killed by $\mathbf{m}$
on the left. The image $\mathcal{X}$
of ${}^1_{0} \mathcal{H}_0^1$ under the
equivalence between ${}^\infty_{\,\,0} \mathcal{H}_0^1$
and $\mathcal{O}_0$
consists of all modules in $\mathcal{O}_0$  killed
by $\mathbf{m}$. In particular, $\mathcal{X}$
contains all simples in $\mathcal{O}_0$ as
well as all Verma modules.

The category ${}^1_{0} \mathcal{H}_0^1$
has an auto-equivalence given by swapping
the sides of the bimodule, see \cite{SoAA}.
This induces an auto-equivalence $\Psi$ of 
$\mathcal{X}$ which sends $L_w$ to $L_{w^{-1}}$,
for any $w\in W$. Note that $L_w\cong L_{w^{-1}}$
if and only if $w$ is an involution. It follows that
$\Psi(P^{\hat{\mathcal{R}}}_e)\cong P^{\hat{\mathcal{R}}}_e$
as well as $\Psi(Q_d)\cong Q_d$, if $d\in S_n$ is an involution
and $\mathcal{R}$ the right cell containing $d$.

\subsection{Relation to Kostant's problem}\label{s3.4}

By adjunction, for $L\in\mathcal{O}_0$, we have 
\begin{equation}\label{eq.adj}
\mathrm{Hom}_{\mathcal{O}}(\theta_x L, L)\cong
\mathrm{Hom}_{\mathcal{H}}(\theta_x,\mathcal{L}(L,L)).
\end{equation}
For $L=L_d$ where $d\in W$ is an involution the bimodule $\mathcal{L}(L,L)$
corresponds, under the above equivalence, to the module $Q_d$,
see \cite[Section~7]{KMM}. The grading on $Q_d$ is such that 
the $L_d$ subquotient is in degree $\mathbf{a}(d)$, which turns
\eqref{eq.adj} into
\begin{displaymath}
\mathrm{hom}_{\mathcal{O}}(\theta_x L_d, L_d \langle i\rangle)\cong
\mathrm{hom}_{\mathcal{H}}(\theta_x,\mathcal{L}(L_d,L_d)\langle i\rangle),
\end{displaymath}
for any $i\in\mathbb{Z}$.
Since $\theta_x$ corresponds to $P_x\in {}^\mathbb{Z}\mathcal O_0$, we obtain:
\begin{equation}\label{eq.Qd}
    \dim \mathrm{hom}(\theta_x L_d\langle i\rangle, L_d) = [Q_d: L_x\langle i\rangle].
\end{equation}
Therefore Kostant's problem is directly related to the problem of
description of the tops of the modules of the form
$\theta_x L_d$, where $x\in W$. Note that 
$\theta_x L_d\neq 0$ if and only if $x\leq_R d$.
However, $\theta_x L_d\neq 0$ does not imply
$\mathrm{Hom}_{\mathcal{O}}(\theta_x L_d, L_d)\neq 0$, in general.

\subsection{Further conditions}\label{s3.5}

For $x,y\in W$, following \cite{KiM}, 
we write $\mathbf{KM}(x,y)=\mathtt{true}$ provided that 
$\theta_x L_y$ is either indecomposable or zero.
We will write $\mathbf{KM}(x,y)=\mathtt{false}$ provided that 
$\theta_x L_y$ is non-zero but decomposable.
We denote by $\mathbf{KM}(*,y)$ the conjunction of
$\mathbf{KM}(x,y)$ over all $x\in W$.

For an involution $d\in W$, following \cite{KMM}, 
we write $\mathbf{Kh}(d)=\mathtt{true}$ provided that,
for all $x,y\in W$ such that $x\neq y$
and both $x\leq_R d$ and $y\leq_R d$, the modules
$\theta_x L_d$ and $\theta_y L_d$ are not isomorphic.

The condition $\mathbf{Kh}$ is constant on left cells,
see \cite[Proposition~22]{CM}.
The condition $\mathbf{KM}$ is constant on left
cells with respect to  the second argument, see \cite[Proposition~2]{CMZ}.

\subsection{Issue}\label{s3.6}

Theorem~8.16 in \cite{KMM} asserts that, for an involution
$d\in W$, the validity of $\mathbf{K}(d)$ is equivalent to
the conjunction of $\mathbf{Kh}(d)$ with $\mathbf{KM}(*,d)$.

That the failure of either $\mathbf{Kh}(d)$ or 
$\mathbf{KM}(*,d)$ implies the failure of $\mathbf{K}(d)$
is fairly easy, see \cite[Proposition~8.14]{KMM}. 
Unfortunately, the proof of the fact that 
the validity of both $\mathbf{Kh}(d)$ with $\mathbf{KM}(*,d)$
implies the validity of $\mathbf{K}(d)$ is based on 
\cite[Proposition~8.15]{KMM} whose proof
contains a gap. Namely, the image of the functor $\Theta$
defined in that proof gives 
$\theta_x\otimes_U \mathrm{Res}^{A_d}_{U_d}(A_d)$,
in general, and not directly 
$\theta_x\otimes_U U_d$ as claimed.
This is a gap in the proof which we cannot repair in the general case
at the moment.

\subsection{Known approaches to Kostant positivity}\label{s3.7}

So far, it is unclear whether the validity of
$\mathbf{Kh}(d)$ and $\mathbf{KM}(*,d)$
is sufficient for $\mathbf{K}(d)$.
Let us therefore recall the known approaches
to establish Kostant positivity.

Let $W^P$ be a parabolic subgroup of $W$ and
$w_0^P$ the longest element in $W^P$. 
Let $w_0$ be the longest element of $W$.
By \cite{GJ}, we have $\mathbf{K}(w_0^Pw_0)=\mathtt{true}$.
This implies both $\mathbf{KM}(*,w_0^Pw_0)=\mathtt{true}$
and $\mathbf{Kh}(w_0^Pw_0)=\mathtt{true}$.

The paper \cite{Kh} also proves the following reduction result:
if $W^P$ is a parabolic subgroup of $W$ and $w\in W^P$,
then $\mathbf{K}^P(w)=\mathbf{K}(ww_0^Pw_0)$, where
$\mathbf{K}^P$ is the value of the Kostant positivity 
function on $w$ with respect to
the Lie algebra corresponding to $W^P$.

\section{Kostant's problem for fully commutative permutations}\label{s4}

\subsection{Fully commutative permutations}\label{s4.1}

Recall that $w\in S_n$ is fully commutative if the shape of its 
Robinson-Schensted correspondent has at most two parts.
Alternatively, $w\in S_n$ is fully commutative if it
avoids the pattern $321$. 

If $w\in S_n$ is fully commutative and  its 
Robinson-Schensted correspondent has shape $\lambda=(\lambda_1,\lambda_2)$,
then $\mathbf{a}(w)=\lambda_2$, see \cite{MT}.

\subsection{Temperley-Lieb algebra}\label{s4.2}

Fully commutative elements are best described by 
Temperley-Lieb diagrams, see \cite{MMM} for details.
For a fully commutative $w\in S_n$, we denote the
corresponding Temperley-Lieb diagram by $\mathbf{d}(w)$.

For example, the permutation $21354$
(one-line notation) has reduced expression $s_1s_3$
(where $s_i$ is the elementary transposition $(i,i+1)$).
Its Robinson-Schensted correspondent  
\begin{displaymath}
\begin{ytableau}
1&3&4\\2&5 
\end{ytableau}
\end{displaymath}
has shape $(3,2)$ and hence this permutation is fully commutative.
The corresponding Temperley-Lieb diagram is:
\begin{displaymath}
\xymatrix@R=5mm@C=5mm{
\bullet\ar@{-}@/_/[r]&\bullet&\bullet\ar@{-}[d]&
\bullet\ar@{-}@/_/[r]&\bullet\\
\bullet\ar@{-}@/^/[r]&\bullet&\bullet&\bullet\ar@{-}@/^/[r]&\bullet
}
\end{displaymath}
The value of the $\mathbf{a}$-function of $w$ coincides with
the number of cups (or caps) in $\mathbf{d}(w)$.
Thus, the above element $21354$ has $\mathbf{a}$-value $2$.

The Temperley-Lieb algebra $\mathbf{TL}_n$ 
is an $\mathbb{Z}[v,v^{-1}]$-algebra spanned by all Temperley-Lieb
diagrams with the usual multiplication in which 
removal of each closed circle contributes a $v+v^{-1}$ factor.
This algebra is isomorphic to the quotient of 
$\mathbf{H}$ by the ideal spanned by all $\underline{H}_w$,
where $w$ is not fully commutative. Under this isomorphism,
the element $\underline{H}_w$, where $w$ is fully commutative,
is sent to $\mathbf{d}(w)$.

We write the product $\underline{H}_x\underline{H}_y$
by placing $\mathbf{d}(x)$ on top of $\mathbf{d}(y)$.
In other words, left $\mathbf{H}$-modules correspond to 
top $\mathbf{TL}_n$-modules while right 
$\mathbf{H}$-modules correspond to bottom $\mathbf{TL}_n$-modules.
With this convention, the left KL order is given by 
cap containment: for two fully commutative elements of $W$,
we have $x\leq_L y$ if and only if each cap of $x$ is a cap of $y$.
Similarly, the right  KL order is given by 
cup containment: for two fully commutative elements of $W$,
we have $x\leq_L y$ if and only if each cup of $x$ is a cup of $y$.

\subsection{Tops of $\theta_x L_y$}\label{s4.3}

Let $y\in S_n$ be fully commutative. If we take some $x\in S_n$
such that $\theta_x L_y\neq 0$, then $x$ is also fully commutative.
Furthermore, by \cite[Theorem~4.11(iii)]{BS}, for such $x$ and $y$,
the module $\theta_x L_y$ has simple top.

\begin{lemma}\label{lem:topxyfc}
Let $x,y\in S_n$ be fully commutative and such that 
$\theta_x L_y\neq 0$. Then the simple top of 
$\theta_x L_y$ lives in degree $-\mathbf{a}(x)$.
\end{lemma}

\begin{proof}
By \cite[Theorem~2.2]{KMM}, the simple top of the module $\theta_x L_y$
lives in some degree $k\leq -\mathbf{a}(x)$.

On the other hand, any top of $\theta_x L_y$ comes from some
top of $\theta_x \theta_y$. The latter corresponds to
the product  $\underline{H}_x\underline{H}_y$ in
$\mathbf{H}$. The relevant fully commutative part 
can be computed via the product 
$\mathbf{d}(x)\mathbf{d}(y)$ in $\mathbf{TL}_n$.

When multiplying two TL-diagrams, the maximal number of 
circles that can be removed does not exceed the minimum
between the number of caps in $x$ and the number of cups in $y$.
Since $\theta_x L_y\neq 0$, we have $\mathbf{a}(x)\leq \mathbf{a}(y)$.
Therefore our minimum is the number of caps in $x$, that is 
$\mathbf{a}(x)$. 

Each removed circle contributes a factor $v+v^{-1}$, hence that
maximal contributed factors could be $(v+v^{-1})^{\mathbf{a}(x)}$.
This implies that the simple top of $\theta_x L_y$
lives in some degree $k\geq -\mathbf{a}(x)$.
Putting the two inequalities together, 
we obtain that this degree must be exactly $-\mathbf{a}(x)$.
\end{proof}

\begin{corollary}\label{cor:LLLfc}
Let $d\in S_n$ be a fully commutative involution.
Then the module $Q_d$ is multiplicity free and, moreover,
for any fully commutative $w\in S_n$, we have
$[Q_d:L_w\langle i\rangle]\neq 0$ implies $i=-\mathbf{a}(w)$. 
\end{corollary}

\begin{proof}
By \eqref{eq.Qd}, the fact that  $Q_d$ is multiplicity free follows from
the mentioned property that all non-zero 
$\theta_x L_y$, for fully commutative $x$ and $y$, have simple tops, and the degree claim follows from Lemma~\ref{lem:topxyfc}.
\end{proof}

\begin{corollary}\label{cor:kmmfc}
Let $d\in S_n$ be a fully commutative involution.
Then $\mathbf{K}(d)$ is equivalent to
the conjunction of $\mathbf{Kh}(d)$ with $\mathbf{KM}(*,d)$.
\end{corollary}

\begin{proof}
We only need to prove that 
the conjunction of $\mathbf{Kh}(d)$ with $\mathbf{KM}(*,d)$
implies $\mathbf{K}(d)$.
Note that $\mathbf{KM}(*,d)$ is always true in this setup due to 
the simple top property \cite[Theorem~4.11(iii)]{BS}.
So, we need to prove that $\mathbf{Kh}(d)$ implies
$\mathbf{K}(d)$, alternatively, that the negation of 
$\mathbf{K}(d)$ implies the negation of $\mathbf{Kh}(d)$.

So, assume that $\mathbf{K}(d)$ fails and let 
$L_x\langle -\mathbf{a}(x)\rangle$ be a simple subquotient
of $Q_d/P_e^{\hat{\mathcal{R}}}$. Please note that the
degree shift $\langle -\mathbf{a}(x)\rangle$ is guaranteed by Corollary~\ref{cor:LLLfc}.
We have $x\leq_R d$. Let $d'$ be the involution in the 
KL right cell of $x$. Then $d'\leq_R d$ and, consequently, 
\begin{displaymath}
1\leq [P_e^{\hat{\mathcal{R}}}:L_{d'}\langle -\mathbf{a}(d')\rangle]
\leq [P_e:L_{d'}\langle -\mathbf{a}(d')\rangle]=1, 
\end{displaymath}
where the first inequality follows from 
$\theta_{d'} P_e^{\hat{\mathcal{R}}}=P_{d'}^{\hat{\mathcal{R}}}\neq 0$.
In particular, we have $x\neq d'$ due to the multiplicity freeness
of $Q_d$, see Corollary~\ref{cor:LLLfc}.

By Lemma \ref{lem.adjunctioncomputation} and Corollary \ref{cor:LLLfc}, 
the degree zero endomorphism and homomorphism spaces between $\theta_xL_d$ 
and $\theta_{d'}L_d$ are one dimensional. Non-negativity of the grading,
see Subsection~\ref{s2.4},
implies that $\theta_xL_d\cong\theta_{d'}L_d$, which establishes the 
failure of $\mathbf{Kh}(d)$.
\end{proof}

\begin{lemma}\label{lem.adjunctioncomputation}
    Let $d',d\in W$ be involutions $d'\sim_R x$. Then we have
    \[\dim \mathrm{hom}(\theta_xL_d,\theta_{d'} L_d)= \dim \mathrm{hom}(\theta_{d'}L_d,\theta_x L_d) =[Q_d:L_x\langle -{a}\rangle]\]
    and 
    \[\dim \mathrm{hom}(\theta_xL_d,\theta_{x} L_d) =[Q_d:L_{d'}\langle -a\rangle],\]
    where $a= \mathbf{a}(d')=\mathbf{a}(x)$.
\end{lemma}
\begin{proof}
    We consider $\mathrm{hom}(\theta_xL_d,\theta_{d'} L_d)$ which is identified with $\mathrm{hom}(\theta_{d'}\theta_xL_d,L_d)$ by adjunction.
    Since $d',x$ belong to the same right cell, the decomposition of $\theta_{d'}\theta_x$ has $ \theta_x\langle -a \rangle $ with multiplicity one and the rest of the summand is of the form $\theta_y\langle j \rangle$ for $j>-a$ (see e.g., \cite[Lemma 3.2]{KMM}).
    By \cite[Theorem~2.2]{KMM}, 
we have
$\mathrm{hom}(\theta_yL_d\langle j\rangle,L_d)=0$
for such shift $j$, so the claim follows by \eqref{eq.Qd}.

The other three hom spaces are computed in the same way (see the second equation in \cite[Lemma 3.2]{KMM} for the case $\mathrm{hom}(\theta_xL_d,\theta_{x} L_d)$ ).
\end{proof}

\subsection{Consequence}\label{s4.4}

Since Lemma \ref{lem.adjunctioncomputation} does not assume full commutativity, Corollary~\ref{cor:kmmfc} can be generalized as follows.

\begin{corollary}\label{cor:nicecase}
Let $d\in S_n$ be an involution. Assume that, for any $w\in S_n$
such that $w\leq_R d$, the top of  the module $\theta_w L_d$ 
is concentrated in degree $-\mathbf{a}(w)$. Moreover, assume that 
the top of $\theta_w L_d$ is simple if $w$ is an involution.
Then $\mathbf{K}(d)$ is equivalent to
the conjunction of $\mathbf{Kh}(d)$ with $\mathbf{KM}(*,d)$.
\end{corollary}

\begin{proof}
We already know that $\mathbf{K}(d)$ implies both 
$\mathbf{Kh}(d)$ and $\mathbf{KM}(*,d)$. Therefore
let us assume that $\mathbf{K}(d)$ fails.

Let $w$ be such that $L_w$ appears in the top of
$Q_d/P_e^{\hat{\mathcal{R}}}$. Then $w\leq_R d$
and, our assumption that the top of $\theta_wL_d$ is in degree $-\mathbf{a}(w)$ implies that 
$L_w$ appears in $Q_d/P_e^{\hat{\mathcal{R}}}$ 
in degree $\mathbf{a}(w)$. Let $d'$ be the involution
in the right KL cells of $w$. Then $\mathbf{a}(w)=\mathbf{a}(d')$.

Since the top of $\theta_{d'}L_d$ is simple by assumptions,
$L_{d'}$ appears with multiplicity one in  $Q_d$.
As $d'\leq_R d$, the module $L_{d'}$ appears with non-zero multiplicity
in $P_e^{\hat{\mathcal{R}}}$. Hence this multiplicity must be one,
which means that $L_{d'}$ does not appear in $Q_d/P_e^{\hat{\mathcal{R}}}$.
In particular, $w\neq d'$.
The adjunction computation in Corollary~\ref{cor:kmmfc} implies
that the degree zero endomorphisms of $\theta_{d'}L_d$ are just scalars.

Let $m$ be the multiplicity of $L_{w}$ in $Q_d$. 
By Lemma~\ref{lem.adjunctioncomputation} 
the degree zero homomorphism space between 
$\theta_{d'}L_d$ and $\theta_w L_d$, in both directions,
is $m$-dimensional. We also have that the 
degree zero endomorphism space of 
$\theta_{d'}L_d$ is one-dimensional, by the assumption that 
the top of $\theta_{d'}L_d$ is simple.
This implies that for any non-zero homomorphism from 
$\theta_{d'}L_d$ to $\theta_w L_d$, there is a non-zero
homomorphism from $\theta_w L_d$ to $\theta_{d'}L_d$
such that the composition is the identity on $\theta_{d'}L_d$.
If $m>1$, then $\theta_w L_d$ contains $m$ copies of $\theta_{d'}L_d$ and thus $\mathbf{KM}(w,d)$ fails. If $\mathbf{KM}(w,d)$ holds, then we necessarily have $m=1$ and moreover $\theta_{d'}L_d\cong \theta_w L_d$, so $\mathbf{Kh}(d)$ fails.
\end{proof}

In fact, the same argument proves the following slightly sharper version:

\begin{corollary}\label{cor:nicecase2}
Let $d\in S_n$ be an involution. Assume that
\begin{enumerate}[$($a$)$]
\item\label{cor:nicecase2.1} 
For any $x\in S_n$
such that $x\leq_R d$, any simple subquotient
$L_x$ in the top of $Q_d$
appears in degree $\mathbf{a}(w)$.
\item\label{cor:nicecase2.2} 
For any involution $x\in S_n$, the module $L_x$ appears with 
multiplicity at most one in  $Q_d$ in degree $\mathbf{a}(x)$. 
\end{enumerate}
Then $\mathbf{K}(d)$ is equivalent to
the conjunction of $\mathbf{Kh}(d)$ with $\mathbf{KM}(*,d)$.
\end{corollary}

\section{Small shapes}\label{s6}
\subsection{Shape $(n-2,1^2)$}\label{s6.1}

\begin{proposition}\label{prop:shape1}
If $w\in S_n$ is such that $\mathbf{sh}(w)=(n-2,1^2)$,
then $\mathbf{K}(w)$ is equivalent to the conjunction of
$\mathbf{Kh}(w)$ with $\mathbf{KM}(*,w)$.
\end{proposition}

\begin{proof}
Due to the left cell invariance, we may assume that $w=d$ is an involution. 
Note that $\mathbf{a}(d)=3$. We have three different cases to consider.

{\bf Case~1.} Assume $d=s_1s_2s_1$ or $d=s_{n-1}s_{n-2}s_{n-1}$.
Due to the symmetry of the root system, it is enough to consider the
case $d=s_1s_2s_1$. In this case, $d=32145\dots n$
(one-line notation) belongs to the
left cell of $w$ which is of the form $w_0^Pw_0$, for the parabolic
subgroup of $S_n$ generated by $s_1,s_2,\dots,s_{n-3}$. Therefore  both this $w$
and $d$ are positive with respect to $\mathbf{K}$, 
$\mathbf{Kh}$ and $\mathbf{KM}$.

{\bf Case~2.} Assume $d=s_is_{i+1}s_i$, where $1<i<n-2$.
Then both $\mathbf{K}$ and $\mathbf{Kh}$ fail for $d$, see
\cite[Theorem~4.3]{CM2}.

{\bf Case~3.} Assume $d=(i,j)$, for $i<j-2$. Consider the module
$Q_d$. It has simple socle $L_d$ in degree $\mathbf{a}(d)=3$.
In this case, any involution $d'$ satisfying $d'\leq_R d$
is, in fact, the longest element of some parabolic subgroup of $W$.
In particular, the module $\theta_{d'}L_d$ has simple top
$L_d$ which is, moreover, concentrated in degree $-\mathbf{a}(d')$.

Let $w\leq_R d$ be such that $L_w$ is a subquotient  of $Q_d$.
If $w=e$, then $\mathbf{a}(e)=0$ and, for $i\in\mathbb{Z}$, we have
\begin{displaymath}
[Q_d:L_e\langle i\rangle]=
[P_e^{\hat{\mathcal{R}}}:L_e\langle i\rangle]=\delta_{i,0},
\end{displaymath}
for obvious reasons.

If $\mathbf{a}(w)=1$, then $L_w$ can only appear in degrees $1$
or $2$ in $Q_d$. The left descent set of $d$ equals $\{s_i,s_{j-1}\}$.
Since $w\leq_R d$, the left descent set of $w$ must be contained
in $\{s_i,s_{j-1}\}$. If $w$ is an involution, this forces $w$
to be either $s_i$ or $s_{j-1}$. Consequently, if $w$ is not 
necessarily an involution, then the left descent set of $w$ is 
either $s_i$ or $s_{j-1}$. In particular, it does not contain the
left descent set of $d$. From \cite[Fact~3.2]{Wa}, it thus follows that
\begin{displaymath}
\mathrm{Ext}^1_{\mathcal{O}}(L_d,L_w)=0 
\end{displaymath}
and therefore $L_w$ cannot appear in degree $2$ in $Q_d$.
So, it can only appear in degree $1=\mathbf{a}(w)$.

If $\mathbf{a}(w)=2$, then $L_w$ can only appear in degree $2$.
Now the claim of our proposition follows from
Corollary~\ref{cor:nicecase}.
\end{proof}

\subsection{Shape $(n-3,2,1)$}\label{s6.2}

In this section, we prove:

\begin{proposition}\label{prop:shape2}
If $w\in S_n$ is such that $\mathbf{sh}(w)=(n-3,2,1)$,
then $\mathbf{K}(w)$ is equivalent to the conjunction of
$\mathbf{Kh}(w)$ with $\mathbf{KM}(*,w)$.
\end{proposition}

As usual, due to the left cell invariance, we may 
assume that $w=d$ is an involution. 
In order to prove this proposition, we do a case-by-case
analysis following the split of all involutions of this
shape as in \cite[Section~5]{CM2}. Before this, we make some
general observations. The $\mathbf{a}$-value of the shape
$(n-3,2,1)$ is $4$. With respect to the dominance order
(which coincides with the two-sided order), here is the
list of all smaller shapes:
\begin{itemize}
\item the shape $(n-2,1,1)$ of $\mathbf{a}$-value $3$;
\item the shape $(n-3,3)$ of $\mathbf{a}$-value $3$;
\item the shape $(n-2,2)$ of $\mathbf{a}$-value $2$;
\item the shape $(n-1,1)$ of $\mathbf{a}$-value $1$;
\item the shape $(n-1)$ of $\mathbf{a}$-value $0$.
\end{itemize}

We seek to apply Corollary~\ref{cor:nicecase2}.

\begin{lemma}\label{lem.shape2cor6b}
    The condition~\eqref{cor:nicecase2.2} in Corollary~\ref{cor:nicecase2} is satisfied for every $d\in S_n$ of shape $(n-3,2,1)$. 
\end{lemma}
\begin{proof}
Let $d$ be an involution of shape $(n-3,2,1)$.
The only element with the $\mathbf{a}$-value $0$
is the identity and it appears with multiplicity $1$
in $Q_d$ in degree $0$.

All involutions $x$ of $\mathbf{a}$-value $1$ are
longest elements of some parabolic subgroups, hence,
for such $x$, the module $\theta_x L_d$ is either zero or
has simple top and this top is in degree $-1$. 

The involutions of $\mathbf{a}$-value 2 and 3 are covered in Lemma~\ref{lem-12357} and Lemma~\ref{lem12357-2} respectively.
The only relevant involution of $\mathbf{a}$-value $4$ is $d$, and the module $Q_d$ has $L_d$ in the simple socle of degree $4$.
\end{proof}

\begin{lemma}\label{lem-12357}
For any involution $x$ of $\mathbf{a}$-value $2$, either
the module $\theta_x L_d$ is zero or $L_d$ appears
with multiplicity one in the top of  $\theta_x L_d$ in degree $-2$.
\end{lemma}

\begin{proof}
Involutions $x$ of $\mathbf{a}$-value $2$ split into 
the longest elements of some parabolic subgroups
and elements of the form $s_is_{i-1}s_{i+1}s_i$.
For the former, the module $\theta_x L_d$ is either zero or
has simple top and this top is in degree $-2$.
For the latter, we note that the length of the element is
even. Hence $L_d$ appears as a subquotient of 
$\theta_x L_d$ only in even degrees. Therefore, as a top, $L_d$ can only appear
in degree $-2$. Note also that, for $x=s_is_{i-1}s_{i+1}s_i$,  we have
\begin{displaymath}
\theta_x=\theta_{s_{i}}\theta_{s_{i-1}}\theta_{s_{i+1}}\theta_{s_{i}}. 
\end{displaymath}
Hence, by adjunction, 
\begin{equation}\label{eq:hom-5}
\mathrm{hom}(\theta_x L_d,L_d\langle 2\rangle)=
\mathrm{hom}(\theta_{s_{i+1}}\theta_{s_{i}} L_d,
\theta_{s_{i-1}}\theta_{s_{i}}L_d\langle 2\rangle).
\end{equation}
Assume $\theta_{s_{i}} L_d\neq 0$.
Then either $\theta_{s_{i+1}}L_d=0$ or $\theta_{s_{i-1}}L_d=0$
or both, for otherwise $\theta_x L_d$ would have $L_d$ in degree $-4$,
which is not possible. If $\theta_{s_{i+1}}L_d=\theta_{s_{i-1}}L_d=0$, then 
the indecomposable module $\theta_{s_{i+1}}\theta_{s_{i}} L_d$ must have
simple top in degree $-1$ while 
the indecomposable module 
$\theta_{s_{i-1}}\theta_{s_{i}} L_d\langle 2\rangle$ must have
simple socle in degree $-1$ as well, so the dimension of 
\eqref{eq:hom-5} is at most $1$. 

If $\theta_{s_{i+1}}L_d\neq 0$,
then we have $\theta_{s_{i-1}}L_d=0$. Therefore the indecomposable module
$\theta_{s_{i-1}}\theta_{s_{i}}L_d\langle 2\rangle$ has simple
socle, say $L_y$, in degree $-1$. In particular, 
$L_y$ appears with multiplicity $1$ in $\theta_{s_{i}}L_d$,
which means that $\dim\mathrm{Ext}^1(L_d,L_y)=1$. Consequently,
$L_d\langle 2\rangle$ appears with multiplicity at most $1$ in 
$\theta_{s_{i-1}}L_y\langle 2\rangle$.
The module 
$\theta_{s_{i+1}}\theta_{s_{i}} L_d$ has a  top 
constituent $L_d$ in degree $-2$ with multiplicity $1$.
By the above, the hom-space from it to the degree $-2$ part of
$\theta_{s_{i-1}}L_y\langle 2\rangle$ is at most one-dimensional.
Assume $\theta_{s_{i+1}}\theta_{s_{i}} L_d$ maps onto
$L_y\langle 1\rangle$. Then, by adjunction
\begin{equation}\label{eq:hom-7}
0\neq \mathrm{hom}(\theta_{s_{i+1}}\theta_{s_{i}} L_d,L_y\langle 1\rangle)= 
\mathrm{hom}(\theta_{s_{i}} L_d,\theta_{s_{i+1}}L_y\langle 1\rangle). 
\end{equation}
The module $\theta_{s_{i}} L_d$ has simple  top $L_d\langle 1\rangle$.
Since we assume $\theta_{s_{i+1}}L_d\neq 0$, the element $d$ has
$s_{i+1}$ in the right descent set. But this means that 
$L_d\langle 1\rangle$ cannot appear in $\theta_{s_{i+1}}L_y\langle 1\rangle$,
which leads to a contradiction that the right hand side of
\eqref{eq:hom-7} must be zero. Therefore 
$\theta_{s_{i+1}}\theta_{s_{i}} L_d$ does not map onto
$L_y\langle 1\rangle$ and hence the dimension of \eqref{eq:hom-5}
is at most $1$. 

Similarly one shows that 
the dimension of \eqref{eq:hom-5} is at most $1$ 
in the case $\theta_{s_{i+1}}L_d=0$
and $\theta_{s_{i-1}}L_d\neq 0$.
\end{proof}

\begin{lemma}\label{lem12357-2}
For  $x$ of $\mathbf{a}$-value $3$, the corresponding $L_x$
can appear in $Q_d$ only in degree $3$ and, in case $x$ is an involution,
with multiplicity at most $1$.
\end{lemma}

\begin{proof}
That $L_x$ can appear in $Q_d$ only in degree $3$
follows from \cite[Theorem~2.2]{KMM}.
If $x$ is an involution of $\mathbf{a}$-value $3$ such that 
$x\leq d$ in the Bruhat order, then $L_x\langle -3\rangle$
appears with multiplicity at most $1$ in $Q_d$ by 
\cite[Theorem~1.4]{Xi}. The case when $x$ is the longest element
of a parabolic subgroup is also clear by the same arguments 
as above. The involutions of $\mathbf{a}$-value $3$ that are
not longest elements of parabolic subgroups split into four
families:
\begin{itemize}
\item transpositions $(i,j)=s_is_{i+1}\dots s_{j-1}\dots s_i$, for $j>i+2$;
\item elements of the form  $s_is_{i+1}s_{i-1}s_is_j$
where the support of $s_j$ is disjoint from the 
support of $s_is_{i+1}s_{i-1}s_i$;
\item elements of the form $s_{i+1}s_{i+2}s_{i-1}s_{i}s_{i+1}s_{i-2}s_{i-1}$;
\item elements of the form $s_{i}s_{i-1}s_{i+1}s_{i-2}s_{i}s_{i+2}s_{i-1}s_{i+1}s_{i}$.
\end{itemize}
For the three last families, we see that the support of the involved elements
is at most $6$. It is known that the Kazhdan-Lusztig $\mu$-function 
has values $0$ or $1$ for $S_n$ with $n\leq 9$, see \cite{MW}.
Therefore, for all these families, $L_x\langle -3\rangle$
appears with multiplicity at most $1$ in $Q_d$.
In the case of the first family, that is, for the transpositions $x=(i,j)$,
where $j>i+2$, the fact that $\mu(x,d)\in\{0,1\}$ follows from \cite[Theorem~2]{SSV}.
\end{proof}

\begin{lemma}\label{lem12357-3}
If $x$ is of $\mathbf{a}$-value $1$, the module $L_x$ 
cannot appear in degree $3$ of $Q_d$.  
\end{lemma}

\begin{proof}
Any such appearance
must represent an extension with the simple socle $L_d$ of $Q_d$
in degree $4$. Since any element of $\mathbf{a}$-value $1$ is boolean,
that is, is a multiplicity free product of simple reflections, the inequality
$d\preceq x$ is impossible, as $d$ is not boolean. Therefore,
for $\mathrm{Ext}^1(L_x,L_d)$ to be non-zero, we must have
$x\preceq d$. If $x$ and $d$ are not Bruhat neighbours,
then, by \cite[Fact~3.2]{Wa}, in the case $\mathrm{Ext}^1(L_x,L_d)\neq 0$,
both the left and the right descent sets of $d$ should be subsets of
the left and the right descent sets of $x$, respectively.
But the left and the right descent sets of $x$ are singletons.
At the same time, the shape of $d$ has three rows and hence its
left and right descent sets contain at least two elements each.
So, this case is not possible.

The only remaining case is the case when $x$ and $d$ are Bruhat neighbours.
We claim that this cannot happen. Indeed, assume $x$ and $d$ are Bruhat
neighbours. Then $d=xt$, where $t=(i,j)$, where $i<j$, is a transposition with
$\ell(d)=\ell(x)+1$. By \cite[Lemma~2.1.4]{BB},  we have $x(i)<x(j)$
and there is no $i<k<j$ such that $x(i)<x(k)<x(j)$. Since $x$
has shape $(n-1,1)$, it must contain an increasing subsequence of
length $n-1$ by \cite{Sch}. By the above condition, this subsequence
intersects the set $\{x(i),x(i+1),x(i+2),\dots,x(j)\}$ in at most two elements.
After our swap, we can still keep at least one of those elements.
Therefore $d$ contains an increasing subsequence of length at least $n-2$
and thus cannot have shape $(n-3,2,1)$.
\end{proof}

By Lemma~\ref{lem.shape2cor6b}, Lemma~\ref{lem12357-2}, and Lemma~\ref{lem12357-3}, in order to apply Corollary~\ref{cor:nicecase2},
it only remains  to check that
\begin{itemize}
\item $L_x$ does not appear in degree $2$ in $Q_d$ when $x$ is of $\mathbf{a}$-value $1$ ;
\item $L_x$ does not appear in degree $3$ in $Q_d$ when $x$ is of $\mathbf{a}$-value $2$ .
\end{itemize}
We do this case by case (going through all possible involutions $d$ of shape $(n-3,2,1)$, which we divide into several) in subsections~\ref{s6.2.1}-\ref{s6.2.7}. For the analysis we note that
an involution of shape $(n-3,2,1)$ is always even and has at most three
simple reflections in its left (the same as right) descent set, which leads to further reductions:

\begin{itemize}
\item ($\mathbf{a}$-value $1$ case)
The only $x$ we need to check
are of the form $s_is_{i+1}\dots s_j$ or $s_is_{i-1}\dots s_j$,
where both $s_i$ and $s_j$ are in the left descent set of $d$.
There are at least two and at most six such elements.
We ignore all odd $x$ as they cannot appear in degree $2$.
\item($\mathbf{a}$-value $2$ case)
There are at most six
right cells that can be right related to $d$, namely those
containing $s_is_j$, for $j\neq i,i\pm 1$, with both
$s_i$ and $s_j$ in the left descent set of $d$, and those
containing $s_is_{i+1}s_{i-1}s_i$ with $s_i$ 
in the left descent set of $d$. So, here we have at most $30$
elements to care about and we can ignore all even 
$x$ as they cannot appear in degree $3$.
\end{itemize}

Also, all Kostant negativity results in \cite{CM2} are proved by showing
that $\mathbf{Kh}$ fails. Hence we only need to look at Kostant positive cases
in \cite{CM2}.

\subsubsection{Type (1)}\label{s6.2.1}

Let $n\geq 5$ and $3\leq i\leq n-2$. Consider now the involution
$d=w_i:=(i-2,i+1)(i-1,i+2)$, here given as a product of two cycles.
The second and the third rows of the RS-correspondent for $w_i$
are:
\begin{displaymath}
\begin{ytableau}
\resizebox{0.7mm}{!}{$i$}&
\resizebox{3mm}{!}{$i+2$}
\\
\resizebox{3mm}{!}{$i+1$} 
\end{ytableau}
\end{displaymath}
The (left and right) descent set of this involution is
$\{s_{i-1},s_i\}$. This is a Kostant positive case that needs to be checked.

The only $\mathbf{a}$-value $1$ elements to check are
$x=s_{i}s_{i-1}$ and $s_{i-1}s_{i}$. We check the first
one, the second one is similar. We have
$\theta_x=\theta_{s_i}\theta_{s_{i-1}}$, the module
$\theta_{s_i}L_d$ has simple top $L_d$ in degree $-1$
and the module $\theta_{x}L_d$ has simple top $L_d$ in degree $-2$.
The application of the adjunction morphism from $\theta_{s_{i-1}}$
to $\theta_e\langle 1\rangle$ surjects $\theta_{x}L_d$ onto
$\theta_{s_i}L_d\langle 1\rangle$. Therefore $L_{x}$ appears in degree
$2$ of $Q_d$ with multiplicity $1$ and, in fact, it appears already 
as a subquotient of $P_e^{\hat{\mathcal{R}}}$ as the adjunction
morphism connects this $L_{x}$ in degree $2$ 
to $L_{s_{i-1}}$ in degree $1$.

There are only two $\mathbf{a}$-value $2$ involutions
that can be potentially right comparable to $d$, namely,
$s_{i-1}s_{i}s_{i+1}s_{i-1}$ and $s_{i}s_{i+1}s_{i-1}s_{i}$.
The unique element in the intersection of the left 
cell of $s_{i-1}s_{i}s_{i+1}s_{i-1}$ and the right cell
of $s_{i}s_{i+1}s_{i-1}s_{i}$ (or vice-versa) is even
and hence can be ignored.

This completes the check for this type.

\subsubsection{Types (2) and (2*)}\label{s6.2.2}

Let $n\geq 5$ and $2\leq i\leq n-3$. Consider the involution
$d=w_i:=(i-1,i+2)(i+1,i+3)$, given as a product of two cycles.
The second and the third rows of the RS-correspondent for $w_i$
are:
\begin{displaymath}
\begin{ytableau}
\resizebox{0.7mm}{!}{$i$}&
\resizebox{3mm}{!}{$i+3$}
\\
\resizebox{3mm}{!}{$i+2$} 
\end{ytableau}
\end{displaymath}
This is type (2). Type (2*) is obtained by conjugation with $w_0$.
We will deal with type (2). This is Kostant negative unless
$i=2$ or $i=n-3$, so we only consider these two cases.
In the Kostant positive case $i=2$, the Kostant positivity 
is deduced using \cite[Theorem~1.1]{Kh}. This Kostant positivity 
implies both $\mathbf{KM}$ and $\mathbf{Kh}$,
so this case is done.

Now assume $i=n-3$. The (left and right) descent set of this involution is
$\{s_{n-4},s_{n-2}\}$.
The only $\mathbf{a}$-value $1$ elements to check in this case are
$s_{n-2}s_{n-3}s_{n-4}$ and $s_{n-4}s_{n-3}s_{n-2}$. Both are odd
so we can ignore them due to Lemma~\ref{lem12357-3}.

We have six elements of $\mathbf{a}$-value $2$ to check in this case, namely,
the two elements

\resizebox{\textwidth}{!}{
$
s_{n-4}s_{n-5}s_{n-3}s_{n-4}s_{n-2}s_{n-3}s_{n-1}s_{n-2}
\quad\text{ and }\quad
s_{n-2}s_{n-3}s_{n-4}s_{n-5}s_{n-1}s_{n-2}s_{n-3}s_{n-4},
$
}

which are even and hence can be ignored, but also the four elements
\begin{displaymath}
s_{n-4}s_{n-3}s_{n-5}s_{n-2}s_{n-4},\quad
s_{n-2}s_{n-3}s_{n-1}s_{n-2}s_{n-4},
\end{displaymath}
\begin{displaymath}
s_{n-4}s_{n-5}s_{n-2}s_{n-3}s_{n-4}\quad\text{ and }\quad
s_{n-4}s_{n-2}s_{n-3}s_{n-1}s_{n-2}
\end{displaymath}
which are odd and require checking.
Note that our involution $d=w_{n-3}$
has reduced expression $s_{n-4}s_{n-2}s_{n-3}s_{n-1}s_{n-2}s_{n-4}$
with  no $s_{n-5}$ in it. Hence it is not Bruhat
comparable to the two elements above that are shorter and have
$s_{n-5}$ as a factor. So, we are left with
\begin{displaymath}
x:=s_{n-2}s_{n-3}s_{n-1}s_{n-2}s_{n-4}
\quad\text{ and }\quad
y:=s_{n-4}s_{n-2}s_{n-3}s_{n-1}s_{n-2}=x^{-1}.
\end{displaymath}
We have $d=ys_{n-4}=s_{n-4}x$ and therefore
\begin{displaymath}
\dim\mathrm{Ext}^1(L_d,L_x)=\dim\mathrm{Ext}^1(L_d,L_y)=1.
\end{displaymath}
Consequently, both $L_x$ and $L_y$ can potentially 
appear in $Q_d$ in degree $3$ with multiplicity at most $1$.

So, let us trace the tops of the modules. We start with $x$.
Note that 
\begin{displaymath}
\theta_x=\theta_{s_{n-2}}\theta_{s_{n-3}}
\theta_{s_{n-1}}\theta_{s_{n-4}}\theta_{s_{n-2}}
\end{displaymath}
by the main result of \cite{BW}.
Applying to $L_d$ the functor $\theta_{s_{n-2}}$ produces a module
with simple top $L_d$ in degree $-1$. Applying now $\theta_{s_{n-3}}$
kills this top and outputs an indecomposable module. This will have top
$L_{x'}$ for $x'=ds_{n-3}\prec d$ in degree $-1$. Now, applying
$\theta_{s_{n-1}}$ moves this top to degree $-2$ and, since 
$\mathrm{Ext}^1(L_d,L_{x'})\neq 0$, produces an element
$L_d$ in degree $-1$. An application of $\theta_{s_{n-4}}$
kills the degree $-2$ part and puts a copy of $L_d$ in that degree.
Since this is an extreme degree, this $L_d$ is there in a top.
The final application of $\theta_{s_{n-2}}$ moves this $L_d$
to degree $-3$. So, $L_x$ does appear in $Q_d$ in degree $3$
but, using the adjunction morphism, this $L_x$ surjects onto
the degree $2$ top component of
$\theta_{s_{n-2}s_{n-3}s_{n-1}s_{n-4}}L_d$ and hence does not
contribute to the top of $Q_d$.

Applying $\Psi$, it follows that $L_y$ does not
contribute to the top of $Q_d$ either.
This completes the check for this type.

\subsubsection{Types (3) and (3*)}\label{s6.2.3}

Let $n\geq 6$ and $2\leq i\leq n-4$ and $i+2<k<n$. Consider the involution
$d=w_{i,k}:=(i-1,k)(k-1,k+1)$, given as a product of two cycles.
The second and the third rows of the RS-correspondent for $w_{i,k}$
are:
\begin{displaymath}
\begin{ytableau}
\resizebox{0.7mm}{!}{$i$}&
\resizebox{3mm}{!}{$k+1$}
\\
\resizebox{0.8mm}{!}{$k$} 
\end{ytableau}
\end{displaymath}
This is type (3). Type (3*) is obtained by conjugation with $w_0$.
We will deal with type (3). The right descent set of 
$w_{i,k}$ is $\{s_{i-1},s_{k-1}\}$ and the reduced expression is 
\begin{equation}\label{eq:redexp5}
s_{i-1}s_{i}\cdots s_{k-3}s_{k-1}s_{k-2}\cdots s_{i-1}s_{k}s_{k-1}.
\end{equation}

The $\mathbf{a}$-value $1$ elements that we need to check are:
\begin{displaymath}
x:=s_{i-1}s_{i}\dots s_{k-1}\quad\text{ and  }\quad
y:=s_{k-1}s_{k-2}\dots s_{i-1}.
\end{displaymath}
Let us start by tracing the tops of $\theta_x L_d$.
Applying $\theta_{s_{i-1}}$ produces a module with simple top
$L_d$ in degree $-1$. Now, applying $\theta_{s_{i}}$ kills
this top and produces a module with simple top $L_{ds_i}$ 
in degree $-1$ (note that $\ell(ds_i)>\ell(d)$). This continues
all the way through the application of  $\theta_{k-3}$ which
produces a module with simple top $L_{ds_is_{i+1}\dots s_{k-3}}$.
This top is still killed at the next step when we apply 
$\theta_{s_{k-2}}$ so we get a new top $L_{w}$ in degree
$-1$ with $w\neq d$ due to right descent comparison.
Consequently, after the next step, we can have $L_d$
at most in degree $-1$. This checks $x$ and an application of 
$\Psi$ gives the corresponding statement for $y$.

There are two $\mathbf{a}$-value $2$ involutions that are
right smaller than our $d$, namely: $s_{i-1}s_{k-1}$
and $s_{k-1}s_{k-2}s_{k}s_{k-1}$. So, up to the action of $\Psi$,
we only need to check one element in the intersection
of the right cell of one of these two involutions and 
the left cell of the other one, namely:
\begin{displaymath}
u:= s_{k-1}s_{k-2}\cdots s_{i-1}s_{k}s_{k-1}.
\end{displaymath}
Note that this reduced expression is a suffix of \eqref{eq:redexp5}.
Removing $u$ from $d$ leaves $s_{i-1}s_{i}\cdots s_{k-3}$
from which it is easy to see that the Bruhat interval $[u,d]$
us a boolean lattice on $\{s_{i-1},s_{i},\dots, s_{k-3}\}$.
Since $k>i+2$, this latter set is not a singleton, that is,
$u$ and $d$ are not Bruhat neighbours. Since the Bruhat 
interval is a boolean lattice, the corresponding Kazhdan-Lusztig
polynomial is trivial 
(see \cite{BB} or combine \cite{WiDM} with \cite{Mar}) 
and hence $\mu(u,d)=0$. This means that
$L_u$ cannot appear in degree $3$ of the module $Q_d$.

This completes our check for this type.

\subsubsection{Types (4) and (4*)}\label{s6.2.4}

Let $n\geq 5$ and $2\leq i\leq n-3$. Consider the involution
$d=w_i:=(i-1,i)(i+1,i+3)$, given as a product of two cycles.
The second and the third rows of the RS-correspondent for $w_i$
are:
\begin{displaymath}
\begin{ytableau}
\resizebox{0.7mm}{!}{$i$}&
\resizebox{3mm}{!}{$i+2$}
\\
\resizebox{3mm}{!}{$i+3$} 
\end{ytableau}
\end{displaymath}
This is type (4*) and type (4) is obtained by
conjugation with $w_0$.

The involutions in this type are Kostant negative, so we have nothing to check.

\subsubsection{Types (5) and (5*)}\label{s6.2.5}

Let $n\geq 6$ and $2\leq i\leq n-4$ and $i+3<k\leq n$. Consider the involution
$d=w_{i,k}:=(i-1,i)(i+1,k)$, given as a product of two cycles.
The second and the third rows of the RS-correspondent for $w_{i,k}$
are:
\begin{displaymath}
\begin{ytableau}
\resizebox{0.8mm}{!}{$i$}&
\resizebox{3mm}{!}{$i+2$}
\\
\resizebox{0.9mm}{!}{$k$} 
\end{ytableau}
\end{displaymath}
This is type (5*) and type (5) is obtained by
conjugation with $w_0$.

The involutions in this type are Kostant negative, so we have nothing to check.

\subsubsection{Types (6) and (6*)}\label{s6.2.6}

Let $n\geq 6$ and $2\leq i\leq n-4$ and $i+2<j<n$. Consider the involution
$d=w_{i,k}:=(i-1,i)(j-1,j+1)$, given as a product of two cycles.
The second and the third rows of the RS-correspondent for $w_{i,k}$
are:
\begin{displaymath}
\begin{ytableau}
\resizebox{0.8mm}{!}{$i$}&
\resizebox{0.8mm}{!}{$j$}
\\
\resizebox{3.2mm}{!}{$j+1$} 
\end{ytableau}
\end{displaymath}
This is type (6*) and type (6) is obtained by
conjugation with $w_0$.

Unless $j+1=n$, the involution above is Kostant negative,
so we have nothing to check.
In the case $j+1=n$, Kostant positivity is established
using \cite[Theorem~1.1]{Kh}. This Kostant positivity 
implies both $\mathbf{KM}$ and $\mathbf{Kh}$,
so this case is also done.

This completes the check for this type.

\subsubsection{Types (7) and (7*)}\label{s6.2.7}

Let $n\geq 7$ and $2\leq i\leq n-5$, $i+2<j<n$ and $j+1<k\leq n$. Consider the involution
$d=w_{i,k}:=(i-1,i)(j-1,k)$, given as a product of two cycles.
The second and the third rows of the RS-correspondent for $w_{i,k}$
are:
\begin{displaymath}
\begin{ytableau}
\resizebox{0.8mm}{!}{$i$}&
\resizebox{0.8mm}{!}{$j$}
\\
\resizebox{0.8mm}{!}{$k$} 
\end{ytableau}
\end{displaymath}
This is type (7*) and type (7) is obtained by
conjugation with $w_0$.
We will deal with type (7*). The right descent set of 
$w_{i,k}$ is $\{s_{i-1},s_{j-1},s_{k-1}\}$ and the reduced expression is 
\begin{displaymath}
s_{i-1}s_{j-1}s_{j}\cdots s_{k-2}s_{k-1}s_{k-2}\cdots s_{j}s_{j-1}.
\end{displaymath}

There are three involutions with $\mathbf{a}$-value $1$ that are right smaller
than $d$, namely: $s_{i-1}$, $s_{j-1}$ and $s_{k-1}$.
Up to taking inverses (for which one can use $\Psi$), the three elements in the
intersection of their left and right cells are:
\begin{displaymath}
y_1= s_{i-1}s_{i}\cdots s_{j-1},\quad
y_2= s_{j-1}s_{j}\cdots s_{k-1}\quad\text{ and }\quad
y_3= s_{i-1}s_{i}\cdots s_{k-1}.
\end{displaymath}
Note that 
\begin{displaymath}
\theta_{y_1}= \theta_{s_{i-1}}\theta_{s_{i}}\cdots \theta_{s_{j-1}}
\end{displaymath}
and similarly for $y_2$ and $y_3$, see \cite{BW}.
Let us track the step by step tops in $\theta_{y_1}L_d$.
Applying $\theta_{s_{i-1}}$ to $L_d$ gives a module with 
simple top $L_d$ in degree $-1$. Applying $\theta_{s_{i}}$ 
kills this top and, by indecomposability, gives a module with 
simple top $L_{ds_{i}}$ in degree $-1$. Applying $\theta_{s_{i+1}}$ 
kills this top and, by indecomposability, gives a module with 
simple top $L_{ds_{i}s_{i+1}}$ in degree $-1$. This pattern
continues all the way and we obtain that the
module $\theta_{y_1}L_d$ is concentrated
in degrees $0,\pm 1$, so here everything is fine.

Let us now track the step by step tops in $\theta_{y_2}L_d$.
Applying $\theta_{s_{j-1}}$ to $L_d$ gives a module with 
simple top $L_d$ in degree $-1$. Applying $\theta_{s_{j}}$ 
kills this top and, by indecomposability, gives a module with 
simple top $L_{ds_{j-1}}$ in degree $-1$. Applying $\theta_{s_{j+1}}$ 
kills this top and, by indecomposability, gives a module with 
simple top $L_{ds_{j-1}s_{j}}$ in degree $-1$. This pattern
continues all the way until the penultimate step where
the degree $-1$ component is $L_{ds_{j-1}s_{j}\dots s_{k-3}}$.
Applying the last $\theta_{s_{k-1}}$ pushes this top to 
degree $-2$. So, we do not have any $L_d$ in degree $-2$
and everything is again fine.

By tracing the two steps one after the other
(not that one instance of $\theta_{j-1}$ in the junction should be ommited), we obtain that
$L_d$ does not appear in 
the degree $-2$ part of $\theta_{y_3}L_d$. Therefore
the $\mathbf{a}$-value $1$ case is fine.

There are three involutions with $\mathbf{a}$-value $2$ that are right smaller
than $d$, namely: $s_{i-1}s_{j-1}$, $s_{i-1}s_{k-1}$ and $s_{j-1}s_{k-1}$.
Up to taking inverses (for which one can use $\Psi$), the three elements in the
intersection of their left and right cells are:
\begin{displaymath}
x_1= s_{i-1}s_{j-1}s_{j}\cdots s_{k-1},\quad
x_2= s_{i-1}s_{i}\cdots s_{j-1}s_{k-1}
\end{displaymath}
and 
\begin{displaymath}
x_3= s_{j-1}s_{j}\cdots s_{k-1}s_{i-1}s_{i}\cdots s_{j-1}.
\end{displaymath}
Note that $x_2$ and $x_3$ are Bruhat incomparable with $d$,
therefore the corresponding simple cannot extend $L_d$ so they
do not appear in degree $3$ of $Q_d$. The element $x_1$ is a prefix of $d$
and deleting it we see that the interval $[x_1,d]$ is boolean on a set
with at least $2$ elements. As in type (3), this implies
$\mu(x_1,d)=0$, so $L_1$ cannot appear in degree $3$ of $Q_d$ either.
Therefore the case of $\mathbf{a}$-value $2$ is fine.

This completes the check for this type.

\section{Small ranks}\label{s5}
In this section, we verify that  $\mathbf{K}(w)$ is equivalent to $\mathbf{Kh}(w)$ and $\mathbf{KM}(*,w)$ for each element $w\in W$, where $W$ is of type $A_n$ for $n\leq 6$. 

\subsection{Type $A_1$}\label{s5.1}

In type $A_1$, we have $\mathfrak{g}=\mathfrak{sl}_2$
and $W\cong S_2=\{e,s\}$. Here both elements are
involutions and both are of the form
$w_0^Pw_0$. So, both are Kostant positive
as well as $\mathbf{Kh}$ and $\mathbf{KM}$ positive.

\subsection{Type $A_2$}\label{s5.2}

In type $A_2$, we have $\mathfrak{g}=\mathfrak{sl}_3$
and $W\cong S_3=\{e,s,t,st,ts,w_0\}$, where $w_0=sts=tst$. 
Here we have four left cells and all of them contain
elements of the form $w_0^Pw_0$: $e,st,ts,w_0$. 
So, all elements in this type are Kostant positive
as well as $\mathbf{Kh}$ and $\mathbf{KM}$ positive.

\subsection{Type $A_3$}\label{s5.3}

We have five partitions of $4$, namely:
$(4)$, $(3,1)$, $(2,2)$, $(2,1^2)$ and $(1^4)$.

The cells corresponding to the partitions 
$(4)$, $(3,1)$, $(2,2)$ consist of fully commutative
elements. In particular, we have Corollary~\ref{cor:kmmfc}
for all such elements. (Here we have our first Kostant negative involution, namely, the element  $s_1s_3$, which is $2143$ in one-line notation.)

The cell corresponding to $(1^4)$ consists of $w_0$,
which is of the form $w_0^Pw_0$ and hence 
this element is Kostant positive.

The cell corresponding to $(2,1^2)$ contains three
left cells and they all contain elements of the form 
$w_0^Pw_0$. Hence all elements in this two-sided cell
are Kostant positive.



\subsection{Type $A_4$}\label{s5.4}

We have seven partitions of $5$, namely:
$(5)$, $(4,1)$, $(3,2)$, $(3,1^2)$, $(2^2,1)$, $(2,1^3)$ and $(1^5)$.

The two-sided cells $(5)$, $(4,1)$, $(3,2)$ consist of 
fully commutative elements and Corollary~\ref{cor:kmmfc} applies. 
The two-sided cells 
$(2,1^3)$ and $(1^5)$ consist of Kostant positive elements
as each of their left cells contain some element of the form $w_0^Pw_0$.
The two-sided cells $(3,1^2)$ and $(2^2,1)$ are is covered by Proposition~\ref{prop:shape1} and Proposition~\ref{prop:shape2} respectively.

\subsection{Type $A_5$}\label{s5.5}

In type $A_5$, we have $\mathfrak{g}=\mathfrak{sl}_6$
and $W\cong S_6$ with simple reflections $s_1,s_2,s_3,s_4,s_5$.
The group $S_6$ has $76$ involutions. Here are the eleven 
partitions of $6$, with the corresponding  $\mathbf{a}$-values
and the number of involutions of the corresponding shape.
Additionally, we record the number of Kostant positive
involutions from \cite[Subsection~10.1]{KMM} 
as well as the number of
left cells hit by the elements of the form
$w_0^Pw_0$ or $sw_0^Pw_0$, where $s$ is a simple
reflection in a proper parabolic subgroup $W^P$
(potentially applied recursively).

\resizebox{\textwidth}{!}{
$
\begin{array}{c||c|c|c|c|c|c|c|c|c|c|c}
\lambda&(6)&(5,1)&(4,2)&(3^2)&(4,1^2)&(3,2,1)&(3,1^3)
&(2^3)&(2^2,1^2)&(2,1^4)&(1^6)\\
\hline\hline
\mathbf{a}(\lambda)&0&1&2&3&3&4&6&6&7&10&15\\
\hline
\# \text{ of invol.}&1&5&9&5&10&16&10&5&9&5&1\\
\hline
\mathbf{K}+&1&5&6&1&8&10&7&1&6&5&1\\
\hline
(s)w_0^Pw_0&0&2&2&1&7&8&7&1&6&5&1
\end{array}
$
}

The fully commutative cases, namely, $(6)$, $(5,1)$, $(4,2)$ and 
$(3^2)$ are dealt with using Section~\ref{s4}.
The shape $(4,1^2)$ is dealt with using Proposition~\ref{prop:shape1}.
The shape $(3,2,1)$ is dealt with using Proposition~\ref{prop:shape2}.
For the shapes $(3,1^3)$, $(2^3)$, $(2^2,1^2)$, $(2,1^4)$, $(1^6)$
we see that the number of claimed Kostant positive elements coincide with 
the number of left cells which can be dealt with using
the approach described in Subsection~\ref{s3.7}.

This verifies the equivalence between $\mathbf{K}$ and the conjunction of $\mathbf{Kh}$ with $\mathbf{KM}$ in type $A_5$ and at the same time confirms \cite[Corollary~10.1]{KMM}.

\subsection{Type $A_6$}\label{s5.6}

In type $A_6$, we have $\mathfrak{g}=\mathfrak{sl}_7$
and $W\cong S_7$ with simple reflections $s_1,s_2,s_3,s_4,s_5,s_6$.
This case is considered in \cite{CM} in detail.
Most of Kostant positive involutions were identified
\begin{itemize}
\item using parabolic induction (i.e. Subsection~\ref{s3.7});
\item using the description in the fully commutative case.
\end{itemize}
In both of these cases, we know that $\mathbf{K}$
is equivalent to the conjunction of $\mathbf{Kh}$
and $\mathbf{KM}$. The remaining involutions are of shape $(5,1,1)$ or $(4,2,1)$, thus are dealt with using Subsections~\ref{s6.1} and~\ref{s6.2}.

\vspace{2mm}

\noindent
H.~K.: Department of Mathematics, Uppsala University, Box. 480,
SE-75106, Uppsala, SWEDEN, email: {\tt hankyung.ko\symbol{64}math.uu.se}

\noindent
V.~M.: Department of Mathematics, Uppsala University, Box. 480,
SE-75106, Uppsala, SWEDEN, email: {\tt mazor\symbol{64}math.uu.se}

%


\begin{thebibliography}{99999999}
\bibitem[BG80]{BG} Bernstein, I. N.; Gelfand, S. I. Tensor products of finite- 
and infinite-dimensional representations of semisimple Lie algebras. Compositio
Math. {\bf 41} (1980), no. 2, 245--285. 
\bibitem[BGG76]{BGG} Bernstein, I. N.; Gelfand, I. M.; Gelfand, S. I. 
A certain category of $\mathfrak{g}$-modules. (Russian) Funkcional. Anal. i 
Prilozen. {\bf 10} (1976), no. 2, 1--8.
\bibitem[BW01]{BW}  Billey, S.; Warrington, G. Kazhdan-Lusztig polynomials for 
321-hexagon-avoiding permutations. J. Algebraic Combin. {\bf 13} (2001), no. 2, 
111--136. 
\bibitem[BB05]{BB}
Bj{\"o}rner, A.; Brenti, F.
Combinatorics of Coxeter groups.
Grad. Texts in Math., {\bf 231},
Springer, New York, 2005. xiv+363 pp.
\bibitem[B-W22]{WiDM}
Blundell, C.; Buesing, L.; Davies, A.; Veli{\v c}kovi{\'c}, P.; Williamson, G.
Towards combinatorial invariance for Kazhdan-Lusztig polynomials
Represent. Theory {\bf 26} (2022), 1145--1191.
\bibitem[BS10]{BS}
Brundan, J.; Stroppel, C.
Highest weight categories arising from Khovanov's diagram algebra. II. Koszulity.
Transform. Groups {\bf 15} (2010), no. 1, 1--45.
\bibitem[CMZ19]{CMZ} Coulembier, K.; Mazorchuk, V.; Zhang, X. Indecomposable 
manipulations with simple modules in category $\mathcal{O}$. Math. Res. Lett. 
{\bf 26} (2019), no. 2, 447--499.
\bibitem[CM25a]{CM3}
Creedon, S.; Mazorchuk, V.
Almost all permutations and involutions are Kostant negative.
Math. Pannon. (N. S.) {\bf 31} (2025), no. 1, 54--59.
\bibitem[CM25b]{CM}
Creedon, S.; Mazorchuk, V.
Hecke combinatorics, K{\aa}hrstr{\"o}m's conditions and Kostant's problem.
Preprint arXiv:2510.05817. To appear in Alg. Combin.
\bibitem[CM26a]{CM2}
Creedon, S.; Mazorchuk, V.
Kostant's problem for permutations of shape $(n-2,1,1)$ and 
$(n-3,2,1)$. Preprint arXiv:2601.19537.
\bibitem[CM26b]{CM4}
Creedon, S.; Mazorchuk, V.
Consecutive patterns, Kostant's problem and type $A_6$.
Internat. J. Algebra Comput. {\bf 36} (2026), no. 4, 317--356.
\bibitem[GJ81]{GJ} Gabber, O.; Joseph, A. On the Bernstein-Gelfand-Gelfand 
resolution and the Duflo sum formula. Compositio Math. {\bf 43} (1981), 
no. 1, 107--131. 
\bibitem[Ge06]{Ge}
Geck, M.
Kazhdan-Lusztig cells and the Murphy basis.
Proc. London Math. Soc. (3) {\bf 93} (2006), no. 3, 635--665.
\bibitem[Hu08]{Hu} Humphreys, J. E. Representations of semisimple Lie algebras in 
the BGG category $\mathcal{O}$. Graduate Studies in Mathematics, {\bf 94}. 
American Mathematical Society, Providence, RI, 2008. xvi+289 pp.
\bibitem[Ja83]{Ja} Jantzen, J. C. Einh{\"u}llende Algebren halbeinfacher 
Lie-Algebren. Ergebnisse der Mathematik und ihrer Grenzgebiete (3), 
{\bf 3}. Springer-Verlag, Berlin, 1983. ii+298 pp. 
\bibitem[Jo80]{Jo} Joseph, A. Kostant's problem, Goldie rank and the 
Gel'fand-Kirillov conjecture. Invent. Math. {\bf 56} (1980), no. 3,  191--213.
\bibitem[KL79]{KL} Kazhdan, D.; Lusztig, G. Representations of Coxeter groups 
and Hecke algebras. Invent. Math. {\bf 53} (1979), no. 2, 165--184.
\bibitem[KiM16]{KiM} Kildetoft, T.; Mazorchuk, V. Parabolic projective 
functors in type $A$. Adv. Math. {\bf 301} (2016), 785--803.
\bibitem[KMM23]{KMM}
Ko, H.; Mazorchuk, V.; Mr{\dj}en, R.
Some homological properties of category $\mathcal{O}$, V.
Int. Math. Res. Not. IMRN {\bf 2023}, no. 4, 3329--3373.
\bibitem[Kh10]{Kh} K{\aa}hrstr{\"o}m, J. Kostant's problem and parabolic 
subgroups. Glasg. Math. J. {\bf 52} (2010), no. 1, 19--32.
\bibitem[KhM10]{KhM} K{\aa}hrstr{\"o}m, J.; Mazorchuk, V. A new approach to 
Kostant's problem. Algebra Number Theory {\bf 4} (2010), no. 3, 231--254.
\bibitem[MMM24]{MMM}
Mackaay, M.; Mazorchuk, V.; Miemietz, V.
Kostant's problem for fully commutative permutations.
Rev. Mat. Iberoam. {\bf 40} (2024), no. 2, 537--563.
\bibitem[Mar06]{Mar}
Marietti, M.
Boolean elements in Kazhdan-Lusztig theory.
J. Algebra {\bf 295} (2006), no. 1, 1--26.
\bibitem[Ma10]{Ma2} Mazorchuk, V. Some homological properties of the 
category $\mathcal{O}$. II. Represent. Theory {\bf 14} (2010), 249--263.
\bibitem[MM11]{MM1} Mazorchuk, V.; Miemietz, V. Cell 2-representations of 
finitary 2-categories. Compos. Math. {\bf 147} (2011), no. 5, 1519--1545.
\bibitem[MS08]{MS} Mazorchuk, V.; Stroppel, C. Categorification of (induced) 
cell modules and the rough structure of generalised Verma modules. 
Adv. Math. {\bf 219} (2008), no. 4, 1363--1426.
\bibitem[MT22]{MT} 
Mazorchuk, V.; Tenner, B.
Intersecting principal Bruhat ideals and grades of simple modules.
Comb. Theory {\bf 2} (2022), no. 1, Paper No. 14, 31 pp.
\bibitem[MW03]{MW}
McLarnan, T.; Warrington, G.
Counterexamples to the $0-1$ conjecture.
Represent. Theory {\bf 7} (2003), 181-195.
\bibitem[Sa01]{Sa} 
Sagan, B. The symmetric group
Grad. Texts in Math., {\bf 203}.
Springer-Verlag, New York, 2001, xvi+238 pp.
\bibitem[Sch61]{Sch}
Schensted, C.
Longest increasing and decreasing subsequences.
Canadian J. Math. {\bf 13} (1961), 179--191.
\bibitem[SSV98]{SSV}
Shapiro, B.; Shapiro, M.; Vainshtein, A.
Kazhdan-Lusztig polynomials for certain varieties of incomplete flags
Discrete Math. {\bf 180} (1998), no. 1-3, 345-355.
\bibitem[So86]{SoAA} Soergel, W.
{\'E}quivalences de certaines cat{\'e}gories de $\mathfrak{g}$-modules.
C. R. Acad. Sci. Paris S{\'e}r. I Math. {\bf 303} (1986), no. 15, 725--728.
\bibitem[So90]{So} Soergel, W. Kategorie $\mathcal{O}$, perverse Garben 
und Moduln {\"u}ber den Koinvarianten zur Weylgruppe. J. Amer. Math. Soc. 
{\bf 3} (1990), no. 2, 421--445. 
\bibitem[So07]{So2}
Soergel, W.
Kazhdan-Lusztig-Polynome und unzerlegbare Bimoduln über Polynomringen.
J. Inst. Math. Jussieu {\bf 6} (2007), no. 3, 501--525.
\bibitem[Wa11]{Wa}
Warrington, G.
Equivalence classes for the $\mu$-coefficient of Kazhdan-Lusztig 
polynomials in $S_n$. Exp. Math. {\bf 20} (2011), no. 4, 457--466.
\bibitem[Xi05]{Xi}
Xi, N.
The leading coefficient of certain Kazhdan-Lusztig 
polynomials of the permutation group $S_n$.
J. Algebra {\bf 285} (2005), no. 1, 136--145.
\end{thebibliography}
\end{document}